\documentclass[11pt]{article}
\usepackage{amssymb,supertabular,latexsym,amsmath,amsthm,graphicx}
\theoremstyle{plain}
\newtheorem{theorem}{Theorem}[section]
\newtheorem{lemma}[theorem]{Lemma}

\theoremstyle{definition}
\newtheorem{definition1}[theorem]{Definition}
\newtheorem{step}{Step}
\newtheorem{case}{Case}

\theoremstyle{definition}

\begin{document}

\title{On the equations defining curves\\ in a polynomial algebra}

\author{Ze Min Zeng\\
 Department of Mathematics,\\
 Washington University in St. Louis,\\
St. Louis, Missouri, 63130\\
Email: \texttt{zmzeng@math.wustl.edu} }

\date{}

\maketitle

\begin{abstract}
Let $A$ be a commutative Noetherian ring of dimension $n$ ($n \ge 3$). Let $I$ be a local complete intersection ideal in $A[T]$ of height $n$. Suppose $I/{I^2}$ is free ${A[T]}/I$-module of rank $n$ and $({A[T]}/I)$ is torsion in $K_0(A[T])$. It is proved in this paper that $I$ is a set theoretic complete intersection ideal in $A[T]$ if one of the following conditions holds: (1) $n$ $\ge 5$, odd; (2) $n$ is even, and $A$ contains the field of rational numbers; (3) $n = 3$ , and $A$ contains the field of rational numbers.
\end{abstract}



\par

\section{Introduction}
Let $A$ be a commutative Noetherian ring of dimension $n$. Let $I$ be a local complete intersection ideal in $A[T]$ of height $n$. In [\ref {Ma2}], Mandal proved the following interesting theorem:
\begin{theorem}
Let $A$ be a commutative Noetherian ring of dimension $n$. Let $I$ be a local complete intersection ideal in $A[T]$ of height $n$. Suppose $I$ contains a monic polynomial. Then $I$ is a set theoretic complete intersection in $A[T]$.
\end{theorem}

Let $A$ be a commutative Noetherian ring of dimension $n$. Let $J$ be a local complete intersection ideal in $A[T]$ of height $n$. By the well known Ferrand-Szpiro construction[\ref {Sz}], there exists a local complete intersection ideal $I$ which is contained in $J$, such that $\sqrt{I} = \sqrt{J}$ and $I/{I^2}$ is free ${A[T]}/I$-module of rank $n$. So to show $J$ is a set theoretic complete intersection in $A[T]$, it suffices to show $I$ is a set theoretic complete intersection in $A[T]$. If $A$ contains the field of rational numbers, in this paper, we are able to find a new sufficient condition for $I$ to be a set theoretic complete intersection, which can be read as follows:
\begin{theorem}
Let $A$ be a commutative Noetherian ring of dimension $n$ ($n \ge 3$) containing the field of rational numbers. Let $I$ be a local complete intersection ideal in $A[T]$ of height $n$. Suppose $I/{I^2}$ is free ${A[T]}/I$-module of rank $n$ and $({A[T]}/I)$ is torsion in $K_0(A[T])$. Then $I$ is a set theoretic complete intersection ideal in $A[T]$.
\end{theorem}

If the dimension of $A$ is odd and bigger than 5, the above theorem also holds without the condition $A$ contains the field of rational numbers. 

All rings in this paper are assumed to be commutative and Noetherian. 
All modules considered are assumed to be finitely generated.

\section{Some preliminary results}
In this section, we collect some well known results that will be used in the third section.

First let us begin with a standard definition.
\begin{definition1} 
Let $A$ be a commutative Noetherian ring. $K_0(A)$ is called the {\it Grothendieck
Group } of $A$, which is defined  by taking the free abelian group generated 
by all finitely generated $A$-modules with finite projective dimension 
modulo the relation generated by 
 $(M_2)=(M_1)+(M_3)$, 
whenever we have an exact sequence of finitely generated $A$-modules of finite projective dimension, 
\[0\to M_1\to M_2\to M_3\to 0.\]
\end{definition1}

The following theorem is due to Boratynski[\ref {Bo}] and Murthy[\ref {Mu}]: 
\begin {theorem}
Let $A$ be a commutative Noetherian ring of dimension $n$, and $I \subseteq A$
be a local complete intersection of height $r$ ($ r \le n $). Suppose
 $I/I^2$ is $A/I$-free with base $\bar{f_1},\dots, \bar{f_r}$, $f_i
\in I$, $\bar {f_i}$ is the class of $f_i$ in $I/I^2$. Let 
$J=I^{(r-1)!} + (f_1, \dots, f_{r-1})$. Then there exists a surjection 
$P \to J $ with $P$ a projective $A$-module of rank $r$, such that 
$(P)- (A^r)= -(A/I) $ in $K_0(A)$. 

\end{theorem}

The next useful lemma which follows is due to Mohan Kumar[\ref {Mk}, Lemma 1].

\begin{lemma}
Let $A$ be a ring and $J \subset A$ a finitely generated ideal. Suppose that $J/J^2$ is generated by $n$ elements. Then for any $a \in A$, the ideal $(J, a)$ is generated by $n+1$ elements.
\end{lemma}

Next we state the cancellation theorem of Plumstead[\ref {Pl}, Theorem 1]:

\begin{theorem}
Let $A$ be a commutative Noetherian ring of dimension $n$, let $P$ and $P'$ be projective $A[T]$-modules with rank $\ge {n+1}$. If $P \oplus A[T] \approx P' \oplus A[T] $, then $P \approx P'$.
\end{theorem}

The following theorem is a result of Bass[\ref {Ba}].

\begin{theorem}
Let $B$ be a ring with dimension $n+1$, $P$ a stably free projective $B$-module of rank $n$, such that $P \oplus B \approx B^{n+1}$.
If $n$ is odd, then $P$ has a unimodular element. 
\end{theorem}

The following theorem is a consequence of a result of Rao[\ref {Ra}, Corollary 2.5] and Quillen's local-global principle[\ref {Qu}, Theorem 1].
\begin{theorem}
Let $A$ be a Noetherian ring of dimension $n$. Suppose $n!$ is invertible in $A$. Then any projective module given by a unimodular row over $A[T]$ of length $n+1$ is extended from $A$. In other words, all stably free $A[T]$-module of rank $n$ are extended from $A$.
\end{theorem}

Let $A$ be a commutative Noetherian ring of dimension $n$ which contains the field of rational
numbers. The Euler class group $E(A)$ of $A$ is defined by Bhatwadekar and Sridharan in [\ref {Bs2}]. Many 
important and interesting theorems are proved in their paper. Let us quote one of them
[\ref {Bs2}, Theorem 4.2]:

\begin{theorem}
Let $A$ be a Noetherian ring of dimension $n \ge 2$ which contains the field of rational
numbers. Let $J$ be an ideal of height $n$ such that $J/J^2$ is generated by $n$ elements, and 
let $\omega_J :(A/J)^n \to J/J^2$ be a local orientation of $J$. Suppose that the image of
 $(J, \omega_J)$ is zero in the Euler class group $E(A)$ of $A$. Then $\omega_J$ is a global 
orientation of $J$.  In other words, $\omega_J$ can be lifted to a surjection $\psi : A^n \to J $.
\end{theorem}  

The following theorem is due to Das[\ref {Da}, Theorem 3.10]
\begin{theorem}
Let $A$ be a Noetherian ring of dimension $n \ge 3 $, containing the field of rational numbers. Let $I \subset A[T]$ be an ideal of height $n$. Suppose there exists a surjection $\phi : A[T]^n \to I/{I^2T}$, such that $\phi \otimes A(T)$ can be lifted to a surjection $\phi' : A(T)^n \to IA(T)$. Then there is a surjection $\psi : A[T]^n \to I$ which lifts $\phi$.
\end{theorem}

\par

\section{Set theoretic complete intersection for curves in $A[T]$}
Let $A$ be a commutative Noetherian ring of dimension $n$. Let $I$ be a local complete intersection ideal of height $n$ in $A[T]$. Assume $I/{I^2}$ is free ${A[T]}/I$-module of rank $n$. If $({A[T]}/I)$ is torsion in $K_0(A[T])$, then by the standard ``thickening'' method, we will show that we can find a  local complete intersection ideal $K$ in $A[T]$, such that $\sqrt{K} = \sqrt{I} $ and $K$  is the surjective image of 
a stably free projective $A[T]$-module $\tilde{P}$ of rank $n$. More precisely:

\begin{lemma}
Let $A$ be a commutative Noetherian ring of dimension $n$ ($n\ge 2$). Let $I$ be a local complete intersection ideal of height $n$ in $A[T]$, such that $I/{I^2}$ is free ${A[T]}/I$-module of rank $n$. If $({A[T]}/I)$ is torsion in $K_0(A[T])$, then there exist a local complete intersection ideal $K$ in $A[T]$, and a surjection $\xi : \tilde{P} \to K $ such that:
\begin{enumerate}
\item $K \subset I$, $\sqrt{K}= \sqrt{I} $ , 
\item $\tilde{P} \oplus A[T] \approx {A[T]}^{n+1}$, 
\item $({A[T]}/K)=0$ in $K_0(A[T])$ and $K/{K^2}$ is free ${A[T]}/K$-module of rank $n$.
\end{enumerate}
\end{lemma}

\begin{proof}
We can find a regular sequence  $f_1, \dots, f_n$ in $I$, such that 
$I=(f_1, \dots, f_n) + I^2 $. Suppose $r({A[T]}/I)=0$ in $K_0(A[T])$. Let
$J=I^r + (f_1, \dots, f_{n-1})$, then $({A[T]}/J)=r({A[T]}/I)=0$ in $K_0(A[T])$ and
$J=( f_1, \dots,f_{n-1}, f_n^r) + J^2 $. Let $K=J^{(n-1)!} + (f_1,\dots, f_{n-1})$, then K is a locally complete intersection ideal in $A[T]$, satisfying $ ({A[T]}/K)=0 $ in $K_0(A[T]) $, 
$\sqrt{I}= \sqrt{K}$ and $K/{K^2}$ is free ${A[T]}/K$-module of rank $n$ generated by the image of $ f_1, \dots,f_{n-1}, f_n^{r(n-1)!} $. By theorem 2.2, there exists a projective $A[T]$-module $\tilde{P}$ of rank $n$ and a surjection: $\tilde{P} \to K $, such that $(\tilde{P}) - (A[T]^n) =0 $ in $K_0(A[T])$. Hence $\tilde{P}$ is a stably free $A[T]$-module of rank $n$. By Plumstead cancellation theorem, $\tilde{P} \oplus A[T] \approx {A[T]}^{n+1} $. The proof of the lemma is complete.

\end{proof}

\begin{theorem}
Let $A$ be a commutative Noetherian ring of dimension $n$ ($n\ge 4$). Let $I$ be a local complete intersection ideal of height $n$ in $A[T]$, such that $I/{I^2}$ is free ${A[T]}/I$-module of rank $n$. Suppose $({A[T]}/I)$ is torsion in $K_0(A[T])$  and $n$ is odd. Then $I$ is a set theoretic complete intersection in $A[T]$.
\end{theorem}

\begin{proof}
By Lemma 3.1, there exist a local complete intersection ideal $K$ in $A[T]$, and a surjection $\xi : \tilde{P} \to K $, such that $\sqrt{K} = \sqrt{I} $ and $\tilde{P}$ is a stably free projective $A[T]$-module of rank $n$, and $\tilde{P} \oplus A[T] \approx {A[T]}^{n+1}$ . Since $n$ is odd, by theorem 2.5, $\tilde{P}$ has a unimodular element. So we can write $\tilde{P} = A[T] \oplus \tilde{Q}$ where $\tilde{Q}$ is a stably free $A[T]$-module of rank $n-1$. Let $K_1= \xi(\tilde{Q})$. By some suitable elementary transformation on $\tilde{P}$, we may assume ht$K_1=n-1$. Let $\xi((1,0)) = x $, then $(x, K_1)=K$. Since $n \ge 4$, by Bass cancellation theorem it is easy to see that $\tilde{Q}/{K_1\tilde{Q}}$ is free $A[T]/K_1$-module of rank $n-1$. Therefore, by lemma 2.3, $K=(K_1, x)$ is $n$ generated, and hence $K$ is a complete intersection. Thus $I$ is  a set theoretic complete intersection in $A[T]$.  
\end{proof}

We need the following lemma to prove our next theorem.  

\begin{lemma} \label {lm4}
Let $B$ be a commutative Noetherian ring of dimension $n+1$ containing a field $k$. Let $I$ be an ideal of height $n$ which is a local complete intersection in $B$, such that $I/{I^2}$ is free $B/I$-module of rank $n$.  Then there exists a regular sequence $f_1, \dots, f_n $ in $B$ and $s_1 \in I^2 $ such that
\begin{enumerate} 
\item $I=(f_1, \dots, f_n, s_1)$, $s_1(1-s_1) \in (f_1, \dots, f_n)$,  $I= (f_1, \dots, f_n) + I^2$, and
\item $ \{f_1, \dots, f_{n-1}, f_n - s_1^2 \}$  is a regular sequence in $B$.
\end{enumerate}
\end{lemma}

\begin{proof}
There exists a regular sequence $f_1, \dots, f_n $ in $B$ such that $I=(f_1, \dots, f_n)+ I^2$. By Nakayama's lemma, there exists $s \in I $ such that  $s(1-s) \in (f_1, \dots, f_n)$ and $I= (f_1, \dots, f_n, s)$. Since $s(1-s) \in (f_1, \dots, f_n)$, we may further assume that $s \in I^2 $. Notice that we can change $s$ by $\prod_{i=1}^{m} (s-b_if_n)$ for any positive integer $m$ and $b_i \in B$. If $p_1, \dots, p_t $ are the maximal elements in Ass$(B/{(f_1, \dots, f_{n-1})})$, then $f_n \notin p_1, \dots, p_t$. 

If $s \in p_1, \dots, p_t$, then $f_n -s^2 \notin p_1, \dots, p_t$, and we are through.

If, say for example, $s \notin p_1$, but $s+bf_n \in p_1$ for some $b \in I$, we replace $s$ by $s(s+bf_n)$ and assume $s \in p_1$. Repeating this procedure (that is, replacing $s$ by $\prod_{i=1}^{m} (s-b_if_n)$)  and reordering $p_i$ where $i \in \{1, \dots, t \}$ if necessary, we may assume that $s \in p_1, \dots, p_r $, $s-bf_n \notin p_k$ for $k >r$ and any $b \in I$.

Since $s \in p_1, \dots, p_r $, $f_n -s^2 \notin  p_1, \dots, p_r$. If $f_n -s^2 \notin  p_{r+1}, \dots, p_t$, then we are done. So by reordering $p_{r+1}, \dots, p_t$, we may assume $f_n -s^2 \notin p_{r+1}, \dots, p_{r+l}$ and $f_n -s^2 \in p_{r+l+1}, \dots, p_t$. Let $\lambda \in I \cap (\cap_{i=1}^{r+l} p_i ) \setminus \cup_{j=r+l+1}^t p_j $ (such $\lambda $ does exist), and $s_1= s+\lambda f_n$.  Then $f_n - s_1^2 = f_n - s^2 - \lambda f_n(2s + \lambda f_n)$, and $f_n - s_1^2 \notin p_1, \dots, p_{r+l}$ by our choice of $\lambda$. 

Now we claim that $f_n - s_1^2 \notin p_{r+l+1}, \dots, p_t$. If $f_n - s_1^2 \in p_j$ for some $j \in \{ r+l+1, \dots, t \}$, then $2s+\lambda f_n \in p_j$. Notice that since $B$ is a commutative ring containing a field $k$, either 2 is invertible in $B$ or 2 is zero in $B$. If 2 is zero in $B$, then $\lambda f_n \in p_j$, which is impossible.  If 2 is invertible in $B$, then $s+(1/2)\lambda f_n \in p_j$, which contradicts that $s-bf_n \notin p_k$ for $k >r$ and any $b \in I$. So the claim follows.

Therefor $f_n - s_1^2 $  is a nonzero divisor in $B/{(f_1, \dots, f_{n-1})}$. By our choice of $s_1$, we have that $I=(f_1, \dots, f_n, s_1)$, $s_1(1-s_1) \in (f_1, \dots, f_n)$, $s_1 \in I^2$, $I= (f_1, \dots, f_n) + I^2$, and $ \{f_1, \dots, f_{n-1}, f_n - s_1^2 \}$  is a regular sequence in $I$.
\end{proof}

Throughout the rest of this paper, we will denote $P \otimes A[T]$ by $P[T]$ and use similar obvious notation.

Now we state our theorem for the case when dimension of $A$ $\ge 4$, even. The proof of this theorem is motivated by [\ref {Ma1} , Proposition 6.2, Mandal's online talk draft].

\begin{theorem}\label{tt}
Let $A$ be a commutative Noetherian ring of dimension $n$ ($n\ge 4$, even) containing the field of rational numbers. Let $I$ be a local complete intersection ideal of height $n$ in $A[T]$. Suppose there is a surjection $\xi : P[T] \to I $, where $P$ is a stably free $A$-module of rank $n$. Then $I$ is a set theoretic complete intersection in $A[T]$.  
\end{theorem}

\begin{proof}
We give the proof of the theorem in several steps

\begin{step}
We want to find two local complete intersection ideals in $A[T]$ of height $n$, which are comaximal with $I$.

Since $I/{I^2}$ is free ${A[T]}/I$-module of rank $n$ and $I$ is a local complete intersection, we may write $I= (f_1, \dots, f_n) + I^2$, where $\{ f_1, \dots, f_n \}$ is a regular sequence in $A[T]$. By lemma 3.3, we may assume the image of $ f_n - s^2 $ in $A[T]/{(f_1, \dots, f_{n-1})}$ is a nonzero divisor, $I=(f_1, \dots, f_n, s)$ and $s(1-s) \in (f_1, \dots, f_n)$. Let $K_1 = (f_1, \dots, f_n, 1-s)$, then $K_1 \cap I =(f_1, \dots, f_n)$. Since $\{f_1, \dots, f_{n-1}, f_n -s^2 \}$ is a regular sequence and $I= (f_1, \dots, f_{n-1}, f_n -s^2 ) + I^2$, we can write $(f_1, \dots, f_{n-1}, f_n -s^2 )=I \cap K_2$ for some local complete intersection ideal $K_2$ in $A[T]$, which is comaximal with $I$. We may assume $K_1, K_2$ are ideals of  height $n$, since if one of them equals $A[T]$, then $I$ is a complete intersection. Let $g=f_n - s^2$, then $gA[T] + K_1 =A[T]$, and hence $I, K_1, K_2$ are pairwise comaximal. Let $I^{(2)}=(f_1, \dots, f_{n-1})+I^2$, then $I^{(2)} \cap K_1\cap K_2 = (f_1, \dots, f_{n-1}, gf_n)$. It is clear that $g = -s^2$  is unit modulo $K_1$ and $f_n = s^2$ is unit modulo $K_2$.
\end{step} 

\begin{step}
We may assume $( \ f_1(0), \dots, f_{n-1}(0), g(0)f_n(0) \ )$ is an ideal of height $n$ in $A$.

To see this, we can replace $T$ by $T + \lambda$ for some suitable $\lambda \in {\Bbb Q}$, and assume $ ( f_1(0), \dots, f_{n-1}(0), g(0)f_n(0) )=A $  or  ht$( f_1(0), \dots,  g(0)f_n(0) )$ $ =n $. But if $ ( \ f_1(0), \dots, f_{n-1}(0), g(0)f_n(0) \ )=A $, then $I(0)= A$. Thus $P$ has a umimodular element. It is then clear that $I$ is a complete intersection in $A[T]$ by using lemma 2.3 as in the proof of theorem 3.2. Therefore we may assume $ (\  f_1(0), \dots, f_{n-1}(0), g(0)f_n(0) \ )$ is an  ideal of hight $n$ in $A$.
\end{step} 

\begin{step}
We may assume $(f_1, \dots, f_{n-1}, gf_n)A(T)$ is an ideal of height $n$ in $A(T)$.

If not, then $I^{(2)}$ contains a monic polynomial and hence so does $I$. From the surjection $\xi$, we see $P(T) $ has a unimodular element when localized at such a monic polynomial. It implies $P[T]$ has a unimodular element by [\ref {Bs1}, theorem 3.4]. Hence $I$ is a complete intersection. 
\end{step} 

\begin{step}
We define a local orientation on $I^{(2)}$.

Since $(f_1, \dots, f_{n-1}, f_n) = I \cap K_1 $ and $K_1$  is a local complete intersection ideal of height $n$ and comaximal with $I$ in $A[T]$, we have two natural surjective homomorphisms
\[ 
\begin{array}{ccc}
  (A[T])^n & \stackrel{\alpha}{\to}    &  I \cap K_1  
\\
  
(A[T])^n & \stackrel{\alpha'}{\to}    &  I \cap K_1

\end{array} \] 

defined by $\alpha(e_i)= f_i$ for $i= 1, \dots,n$, and
 $ \alpha'(e_j)=f_j $ for $j=1,\dots, n-1$, $\alpha'(e_n)=-f_n$
where $e_1, \dots, e_n$ is the standard basis of $A^n$ . 
Since $(f_1, \dots, f_{n-1}, g)= I \cap K_2 $ and $I^{(2)} \cap K_1\cap K_2 = (f_1, \dots, f_{n-1}, gf_n)$, where $I, K_1, K_2$ are pairwise comaximal, we also have two natural  surjective homomorphisms

\[ 
\begin{array}{ccc}
  (A[T])^n & \stackrel{\beta}{\to}    &  I \cap K_2  
\\
  
({A[T]})^n & \stackrel{\gamma}{\to}    &  I^{(2)} \cap K_1 \cap K_2

\end{array} \]

defined by $\beta(e_i)= f_i$ for $i= 1, \dots,n-1$, $\beta(e_n)=g$, and
 $ \gamma(e_j)=f_j $ for $j=1,\dots, n-1$, $\gamma(e_n)=gf_n$. 
We define a local orientation $\omega$ of $I^{(2)}$ as $\omega=\gamma\otimes A[T]/{I^{(2)}}$.
\end{step} 

\begin{step}
We want to show $\omega \otimes A[T]/T$ can be lifted to a surjection $\theta : A^n \to I^{(2)}(0)$.

Let $\xi(0)= \xi \otimes A[T]/T$, $\omega(0)= \omega \otimes A[T]/T$, $\omega_I(0)= \alpha \otimes A[T]/I \otimes A[T]/T$, $\omega_{K_1}(0)= \alpha \otimes A[T]/K_1 \otimes A[T]/T $, and $ \omega_{K_2}(0)= \beta  \otimes A[T]/K_2 \otimes A[T]/T$. We show $(I^{(2)}(0), \omega(0)) = 0 $ in $E(A)$, the Euler class group of $A$. There are four cases: \par

\begin{case}
 $K_1(0)$ and $K_2(0)$ are ideals of height $n$ in $A$. \par
From the surjection $\alpha'$, we have the following relation in $E(A)$:
\begin{enumerate} 

 \item     $ (I(0), -\omega_I(0)) + (K_1(0) , -\omega_{K_1}(0))=0 $. \par

 From the surjection $\beta$, we have $(I(0), \omega_I(0)) + (K_2(0) , \omega_{K_2}(0))=0 $ in $E(A)$.  Since $g(0) = -s^2(0)$  is unit modulo $K_1(0)$ and $f_n(0)= s^2(0)$ is unit modulo $K_2(0)$, from the surjections $\beta$ and $\gamma$ and lemma 5.4 in [\ref {Bs2}], we have the following two relations in $E(A)$: 

   \item    $(I(0), \omega_I(0)) + (K_2(0) , \omega_{K_2}(0))=0 $, and 

   \item    $(I^{(2)}(0), \omega(0)) + (K_1(0) , -s^2(0)\omega_{K_1}(0))+ (K_2(0) , s^2(0)\omega_{K_2}(0))=(I^{(2)}(0), \omega(0)) + (K_1(0) , -\omega_{K_1}(0))+ (K_2(0) , \omega_{K_2}(0))=0$. 

\end{enumerate}

From the above three relations, we have  $ (I(0), \omega_I(0)) + (I(0), -\omega_I(0)) = (I^{(2)}(0), \omega(0)) $ in $E(A)$. 

Since $I(0)$ is the image of $\xi(0): P \to I(0) $ and $P$ is a stably free $A$-module, $(I(0), \omega_I(0)) + (I(0), -\omega_I(0))=0 $ in $E(A)$ by [\ref {Bs2}, Proposition 6.2 and Corollary 7.9]. Thus $(I^{(2)}(0), \omega(0)) = 0$ in $E(A)$.  
\end{case}
       
\begin{case}
$K_1(0)=A$ and $K_2(0)$ is an ideal of height $n$ in $A$. \par
From the surjection $\alpha$, we have $(I(0), \omega_I(0))=0 $ in $E(A)$. Hence from the surjection $\beta$, we have $(K_2(0) , \omega_{K_2}(0)) =0 $ in $E(A)$.  Furthermore, the surjection $\gamma$ implies $ 0= (I^{(2)}(0), \omega(0)) + (K_2(0) , s^2(0)\omega_{K_2}(0))= (I^{(2)}(0), \omega(0)) + (K_2(0) , \omega_{K_2}(0))$. So it is easy to see that $(I^{(2)}(0), \omega(0)) = 0 $ in $E(A)$.  
\end{case}

\begin{case} 
$K_2(0)=A$ and $K_1(0)$ is an ideal of height $n$ in $A$. \par
Similar proof as in case (2).
\end{case}

\begin{case}
$K_1(0)=A$ and $K_2(0)=A$. \par
In this case, $(I^{(2)}(0), \omega(0)) = 0 $ in $E(A)$ by definition.
\end{case}

Therefore $\omega \otimes A[T]/T$ can be lifted to a surjection $\theta : A^n \to I^{(2)}(0)$ by theorem 2.7. And thus $\omega $ induces a surjection $ \phi : A[T]^n \to I^{(2)} /{ { I^{(2)} }^2T }$ by [\ref {Bs3}, Remark 3.9].

\end{step}

\begin{step}
We now show $\omega $ can be lifted to a surjection $\psi : (A[T])^n \to I^{(2)} $.
\par
By a similar argument as in step 5 and considering the Euler class group of $E(A(T))$ (Notice that dim$A(T)=n$), it is easy to see $\omega \otimes A(T)$ can be lifted to a surjection $\phi' : A(T) \to I^{(2)}A(T)$. Thus by theorem 2.8, the surjection $ \phi : A[T]^n \to I^{(2)}/{{I^{(2)}}^2T}$ can be lifted to a surjection $ \psi : A[T]^n \to I^{(2)} $.
\end{step}
\par
Therefore, $I^{(2)}$ is a complete intersection, and thus $I$ is a set theoretic complete intersection ideal in $A[T]$. The proof of the theorem is complete.
\end{proof}

Now we are able to give a proof of theorem 1.2 we stated in the Introduction in the case when $n$ is bigger than 4 and even. 

\begin{theorem}
Let $A$ be a commutative Noetherian ring of dimension $n$ ($n\ge 4$, even) containing the field of rational numbers. Let $I$ be a local complete intersection ideal of height $n$ in $A[T]$, such that $I/{I^2}$ is free ${A[T]}/I$-module of rank $n$. Suppose $({A[T]}/I)$ is torsion in $K_0(A[T])$. Then $I$ is a set theoretic complete intersection in $A[T]$.
\end{theorem}

\begin{proof}
By lemma 3.1, there exist a local complete intersection ideal $K$ in $A[T]$ with height $n$, and a surjection $\xi : \tilde{P} \to K $ such that:
\begin{enumerate}
\item $K \subset I$, $\sqrt{K}= \sqrt{I} $ , 
\item $\tilde{P} \oplus A[T] \approx {A[T]}^{n+1}$. 
\end{enumerate}
By theorem 2.6. $\tilde{P}$ is extended from $A$, and thus $\tilde{P} \approx P[T]$, for some  stably free $A$-module  $P$ of rank $n$. By theorem \ref {tt}, $K$ is a set theoretic complete intersection in $A[T]$. Therefore $I$ is a set theoretic complete intersection in $A[T]$. The proof of the theorem is complete.
\end{proof}

Finally, let us state our theorem for the case when dim$A$=3:  
\begin{theorem}
Let $A$ be a commutative Noetherian ring of dimension 3 containing the field of rational numbers. Let $I$ be a local complete intersection ideal of height $3$ in $A[T]$, such that $I/{I^2}$ is free ${A[T]}/I$-module of rank $3$. Suppose $({A[T]}/I)$ is torsion in $K_0(A[T])$. Then $I$ is a set theoretic complete intersection in $A[T]$.
\end{theorem}
\begin{proof}
By lemma 3.1, there exist a local complete intersection ideal $K$ in $A[T]$ with height $3$, and a surjection $\xi : P \to K $ such that:
\begin{enumerate}
\item $K \subset I$, $\sqrt{K}= \sqrt{I} $ , 
\item $ P \oplus A[T] \approx {A[T]}^4$. 
\end{enumerate}
From the surjection $\xi : P \to K $ and any given trivalization $\chi$ of $\bigwedge ^{3} P$, we get an element $(K, \omega_K) \in E(A[T])$ and $e(P, \chi)= (K, \omega_K)$ in $ E(A[T])$ where $(K, \omega_K)$ is obtained from the pair $(\xi, \chi)$. By theorem 2.5, $P$ has a unimodular element. Thus  $e(P, \chi) =0 $ in $ E(A[T])$ by [\ref {Da}, corollary 4.11]. It follows that $(K, \omega_K)=0$ in $ E(A[T])$, and thus $K$ is generated by 3 elements by [\ref {Da}, theorem 4.7]. Therefore $I$ is a set theoretic complete intersection in $A[T]$.

\end{proof}

\section*{Acknowledgement}
I sincerely thank my advisor, Professor N. Mohan Kumar, for his guidance and many valuable suggestions, that have made this paper possible in its present form. I also wish to thank Dr. G.V. Ravindra for his encouragement.

\par

\bibliographystyle{amsplain}

\begin{thebibliography}{10}

\bibitem {Ba} \label{Ba} H. Bass, \textit{Modules which support non-singular forms}, J. Algebra 13 (1969),  246--252.
\bibitem{Bs1} \label{Bs1} S.M. Bhatwadekar, R. Sriharan, \textit{On a question of Roitman}, J. Ramanujan. Math. Society 16 (2001), 45--61.
\bibitem{Bs2} \label{Bs2} S.M. Bhatwadekar, R. Sriharan, \textit{The Euler class group of a Noetherian ring}, Compositio Math. 122 (2000), 183--222.
\bibitem{Bs3} \label{Bs3} S.M. Bhatwadekar, R. Sriharan, \textit{Projective generation of curves in polynomial extension of an affine domain and a question of Nori}, Invent. Math. 133 (1998), 161--1922.
\bibitem{Bo} \label{Bo} M. Boratynski, \textit{A note on set-theoretic complete intersection ideals}, J. Algebra 54 (1978), 1--5.
\bibitem {Da} \label{Da} M. K. Das, \textit{The Euler class group of a polynomial algebra}, J. Algebra 264   (2003), 582--612.
\bibitem{Ma1} \label{Ma1} S. Mandal, \textit{Euler cycles, http://www.math.ku.edu/\~{}mandal/talks/talkEuler.pdf}, (2005)
\bibitem{Ma2} \label{Ma2} S. Mandal, \textit{On set theoretic intersection in affine spaces}, J. Pure App. Algebra 51 (1988), 267--275.
\bibitem{Mk} \label{Mk} Mohan Kumar, N., \textit{Complete intersections}, J. Math. Kyoto Univ. 17 (1977),  533--538.
\bibitem{Mu} \label{Mu} M.P. Murthy, \textit{ Zero cycles and projective modules}, Ann. Math. 140 (1994), 405--434.
\bibitem{Pl} \label{Pl} B. Plumstead, \textit{The conjecture of Eisenbud and Evans}, Amer. J. Math. 105 (1983), 1417--1433.
\bibitem{Qu} \label{Qu} D. Quillen,  \textit{Projective modules over polynomial rings}, Invent. Math. 36(1976), 167--171.
\bibitem{Ra} \label{Ra} R. A. Rao, \textit{The Bass-Quillen conjecture in dimension but characteristic $\ne$ 2, 3 via a question of A. Suslin}, Invent. Math. 93 (1988), 609--618.
\bibitem{Sz} \label{Sz} L. Szpiro, \textit{Equations defining space curves}, Tata Institue, Bombay, 1979.

\end{thebibliography}

\end{document}